\documentclass[a4paper,11pt]{amsart}
\usepackage[latin1]{inputenc}
\usepackage[english]{babel}
\usepackage[T1]{fontenc}
\usepackage{amsmath,amssymb,amsfonts}
\usepackage[all]{xy}

\newcommand{\Z}[0]{\mathbb{Z}}

\newcommand{\C}[0]{\mathbb{C}}
\newcommand{\Q}[0]{\mathbb{Q}}

\newcommand{\B}[1]{\boldsymbol{#1}}

\newcommand{\g}[1]{\mathfrak{#1}}

\newcommand{\sst}[0]{\subset}
\newcommand{\cl}[1]{\mathcal{#1}}

\newcommand{\into}[0]{\hookrightarrow}

\newcommand{\diag}[0]{\Sigma}
\newcommand{\im}[0]{\mbox{im }}

\renewcommand{\qed}[0]{\hspace{\stretch{1}}$\Box$}
\newcommand{\eq}[1][r]
       {\ar@<-3pt>@{-}[#1]
        \ar@<-1pt>@{}[#1]|<{}="gauche"
        \ar@<+0pt>@{}[#1]|-{}="milieu"
        \ar@<+1pt>@{}[#1]|>{}="droite"
        \ar@/^2pt/@{-}"gauche";"milieu"
        \ar@/_2pt/@{-}"milieu";"droite"}
\newcommand{\imm}[1][r] {\ar@{^{(}->}[#1]}

\newtheorem{df}{Definition}[section]
\newtheorem{teo}{Theorem}
\newtheorem{prop}[df]{Proposition}
\newtheorem{lem}[df]{Lemma}

\newtheorem{rem}[df]{Remark}

\newenvironment{dm}[1][Proof.]{\begin{trivlist} \item[\hskip
    \labelsep {\bfseries #1}]}{\end{trivlist}}

\title[On the cohomology of Artin groups in local systems]{On the cohomology of Artin groups in local systems and the associated Milnor fiber}
\author {Filippo Callegaro}
\address{Scuola Normale Superiore, 56100 Pisa, Italy}
\email{f.callegaro@sns.it}
\date {September 2003}

\begin{document}

\begin{abstract}
Let $W$ be a finite irreducible Coxeter group and let $\B{X}_W$ be
the classifying space for $G_W$, the associated Artin group. If
$A$ is a commutative unitary ring, we consider the two local
systems $\cl{L}_q$ and $\cl{L}_q'$ over $\B{X}_W$, respectively
over the modules $A[q,q^{-1}]$ and $A[[q,q^{-1}]]$, given by
sending each standard generator of $G_W$ into the automorphism
given by the multiplication by $q$. We show that
$H^*(\B{X}_W,\cl{L}_q') = H^{*+1}(\B{X}_W,\cl{L}_q)$ and we
generalize this relation to a particular class of algebraic
complexes. We remark that $H^*(\B{X}_W,\cl{L}_q')$
is equal to the cohomology with trivial coefficients $A$ of the Milnor fiber of
the discriminant bundle of the associated reflection group.
\end{abstract}

\maketitle

%\vspace{1cm}

\section*{Introduction}
Let $W$ be a finite irreducible Coxeter group (with Coxeter system 
$(W,S)$) and let $G_W$ be the associated Artin group. 
Recall that if $W = <s,s \in S \mid {(ss')}^{m(s,s')}=e>$ 
is the standard presentation for the Coxeter group, 
then the standard presentation for $G_W$ is given by 
$$
<g_s, s \in S \mid \overbrace{g_sg_{s'}g_sg_{s'}\cdots}^{m(s,s')
\hbox{terms}} = \overbrace{g_{s'}g_sg_{s'}g_s\cdots}^{m(s,s')
\hbox{terms}} \hbox{ for } s \neq s', m(s,s') \neq +\infty> 
$$
(see \cite{br}, \cite{brs} and \cite{deligne}). We call $\B{X}_W$ the
classifying space for $G_W$. Let $A$ be a commutative unitary ring;
we consider a particular local system $\cl{L}_q$ over $\B{X}_W$
with coefficients the ring $A[q,q^{-1}]$, where each standard generator of
$G_W$ acts as $q$-multiplication. Moreover let $\cl{L}_q'$ be the local system which 
is constructed in a similar way over the module $A[[q,q^{-1}]]$.

The cohomology groups $H^*(\B{X}_W,\cl{L}_q')$ have an interesting
geometrical interpretation, in fact they are equal to the
cohomology groups (with trivial coefficients over the ring $A$) of
the Milnor fiber $\B{F}_W$ of the discriminant singularity associated to $W$
(see section \ref{s:appl}). 
From a straightforward application of the Shapiro Lemma (\cite{brown}) it is known that the homology groups $H_*(\B{X}_W,\cl{L}_q)$ are equal to the homology groups of $\B{F}_W$ with coefficients over the ring $A$ (the argument is the same as that used in \cite{cs} for the homology of arrangements of hyperplanes).
%Is is also known that the homology groups
%$H_*(\B{X}_W,\cl{L}_q)$ are equal to the homology groups of $\B{F}_W$ (see \cite{cs}).

The cohomology groups $H^*(\B{X}_W,\cl{L}_q)$ and
$H^*(\B{X}_W,\cl{L}_q')$ can be computed by means of an algebraic
complex described in \cite{salv1}; in this paper we show (see equation (\ref{e:terza})) that
these groups coincide modulo an index shift, that is
$$H^*(\B{X}_W,\cl{L}_q') = H^{*+1}(\B{X}_W,\cl{L}_q). $$
As a consequence we can use $\cl{L}_q$ to compute $H^*(\B{F}_W,A)$. In the
special case when $A = \Q$, and so the ring $A[q,q^{-1}]$ is a
PID, the equality has already been observed (\cite{d}) by Corrado De Concini.
We also give a generalization of this fact, extending the result
to a particular class of algebraic complexes including those
described by Salvetti in \cite{salv1}.

In section \ref{main} we give a precise formulation of the claim
in an algebraic form and we give a proof of it by using spectral
sequences. In section \ref{s:appl} we show how the algebraic result
applies to the cohomology of Artin groups.

\section{Main theorem} \label{main}

\begin{rem} \label{r:indu}
\upshape
Let $(C_1,d)$ be a graduated complex and let $C_3 \sst C_2 \sst C_1$ be inclusions of graduate complexes. Denote by
$d_{ij}:C_i/C_j \to C_i/C_j$ the induced {co\-boun\-da\-ry} on the quotient complex ($1 \leq i < j \leq 3$). There is an obvious exact sequence of complexes:
$$
0 \to C_2/C_3 \into C_1/C_3 \stackrel{\pi}{\rightarrow} C_1/C_2 \to 0
$$

When $d_{12}$ and $d_{23}$ vanish (for example if the complexes are trivial in all degrees except exactly one)
we get that $H^*(C_1/C_2)= C_1/C_2$ and $H^*(C_2/C_3)$ $=$ $C_2/C_3$, so the differential $H^*(C_1/C_2)$ $\to H^*(C_2/C_3)$ of the long exact sequence associated to the above sequence gives a map
$$
\overline{d}: C_1/C_2 \to C_2/C_3.
$$
In the following we call this map \emph{induced differential}.
\end{rem}

Let $A$ be a commutative unitary ring. In this section we indicate by $R=A[q,q^{-1}]$, the ring of Laurent
polynomials with coefficients in $A$ and by $M$ the $R$-module $A[[q,q^{-1}]]$.
Let $(C^*, d^*)$ be a graduate cochain complex, with $C^*$ an $R$-module and $d^*$ an $R$-linear map.
%We recall that if we have a triple of graduate complexes $(C_1, C_2, C_3)$ with $C_1 \sst C_2 \sst C_3$ with 
%coboundary $d$ such that $dC_i \sst C_i$ for all $i$, than the map $d$ induce map $d_{ij}:C_i/C_j \to C_i/C_j$, 
%$1 \leq j < j \leq 3$ which are still coboundary. So we can consider the following short exact sequence of 
%complexes:
%$$
%0 \to C_2/C_3 \into C_1/C_3 \stackrel{\pi}{\rightarrow} C_1/C_2 \to 0
%$$
%In the associated long exact sequence we call "induced differential" the map
%$$
%\overline{d}: H^*(C_2/C_3) \to H^*(C_1/C_2).
%$$
%If we suppose that $d_{12}$ and $d_{23}$ are the null map (for example if each of the complex is trivial, 
%except that in a degree) we get than $H^*(C_1/C_2)= C_1/C_2$ and $H^*(C_2/C_3)= C_2/C_3$ , so we have a map
%$$
%\overline{d}: C_2/C_3 \to C_1/C_2.
%$$
We give the following recursive definition:

\begin{df}
The complex $(C^*, d^*)$ is called \emph{well filtered} if $C^*$ is a free finitely generated $R$-module,
$C^* \neq R$
and moreover, if
%$C^* \neq R$ e
$C^* \neq {0}$, the following conditions are satisfied:

a) $C^*$ is a filtered complex with a decreasing filtration $F$ which is {com\-pa\-ti\-ble} with the coboundary map $d^*$ and
such that $F_0C = C^*$ and $F_{n+1}C = \{ 0 \}$ for an integer $n > 0$ ;

b) $F_nC = (F_nC)^n \simeq F_{n-1}C/F_nC = (F_{n-1}C/F_nC)^{n-1} \simeq R$;

c) the induced differential $\overline{d}: F_{n-1}C/F_nC \to F_nC/F_{n+1}C$ (following from condition (b) and Remark \ref{r:indu}) corresponds to the multiplication by
a non-zero polynomial $p \in R$ with first and last non-zero coefficients invertible in $A$;

d) for all integer $i \neq n-1,n$ the induced complex $((F_iC/F_{i+1}C)^*, d_i^*)$ is a \emph{well filtered} complex.
\end{df}

In the following when we consider a well filtered complex we always suppose to have also a filtration $F$
as above. We write $(C^*_M, d^*_M)$ for the complex $C^* \otimes_R M$ with the natural induced graduation and
coboundary.

\begin{teo} \label{t:centrale}
Let $(C^*,d^*)$ be a well filtered complex. We have the following isomorphism:
$$
H^{*+1}(C^*) \simeq H^{*}(C^*_M).
$$
\end{teo}

In order to proof this fact we need two preliminary lemmas.

As a first step let us consider the natural inclusion of $R$-modules $R \into M$. 
We have the short exact sequence of $R$-modules:
$$
0 \to R \into M \to M' \to 0
$$
where $M' = M/R$. 
We indicate by ${C'}^*$ the complex $C^* \otimes_R M'$ and we consider the complexes
$C^*$, $C^*_M$, ${C'}^* $.
In a similar way we have the following short exact sequence of $R$-modules:
$$
0 \to C^* \stackrel{i}{\into} C_M^* \stackrel{\pi}{\to} {C'}^* \to 0.
$$
Since the maps $i$ and $\pi$ commute with the coboundary maps, we actually have a short exact
sequence of complexes. So we obtain the following long exact sequence:

%\begin{center}
%\begin{tabular}{p{\stretch{1}cp{\stretch{1}}}
\begin{tabular}{rcl}
\xymatrix
@R=1pc @C=1pc
{ \\ (*) \\ \\ }
& 
\xymatrix
@R=1pc @C=1pc
{
& & & \cdots \ar[r]^-{\pi^*}& H^{i-1}({C'}^*) \ar[r]^-{\delta^*} & \\
& \ar[r]^-{\delta^*} & H^i(C^*) \ar[r]^-{i^*} & H^i(C^*_M) \ar[r]^-{\pi^*} & H^i({C'}^*) \ar[r]^-{\delta^*} & \\
& \ar[r]^-{\delta^*} & H^{i+1}(C^*) \ar[r]^-{i^*} & \cdots & & \\ } 
&
\end{tabular}
%\end{center}

\begin{lem} \label{l:iso}
Let $(C^*,d^*)$ be a well filtered complex. With the notation given above we have:
$$
H^i({C'}^*) \simeq H^i(C^*_M) \oplus H^i(C^*_M)
$$
\end{lem}

\begin{dm}
The $R$-module $M'$ splits into the sum of two modules in the following way:
$$
M' = M'_+ \oplus M'_-
$$
where $M'_+ = M/(A[q][[q^{-1}]]), M'_- = M/(A[q^{-1}][[q]])$. In a similar way we get the splitting
$$
{C'}^* = {C'}^*_+ \oplus {C'}^*_-.
$$
Moreover ${C'}^*_+$ and ${C'}^*_-$ are invariant for the coboundary induced by $d^*$,
so the cohomology also splits:
$$
H^*({C'}^*) = H^*({C'}^*_+) \oplus H^*({C'}^*_-).
$$
We want to show that the quotient projection $\pi_+: C_M^* \to {C'}^*_+$ induces an isomorphism
$\pi_{+}^*$ in cohomology. We will prove this by induction on the number of generators of
$C^*$ as a free $R$-module.

%Se $C^* = \{0\}$ o $C^* = R$ il fatto ?ovvio. Supponiamo quindi che $C^*$ abbia $m$ generatori,
%con $m > 1$. I complessi $((F_iC/F_{i+1}C)^*, d_i^*)$ hanno un numero minore di generatori e quindi
%possiamo supporre per ipotesi induttiva che la mappa ${\pi_i}_+$, definita analogamente a $\pi_+$,
%induca un isomorfismo in coomologia per tutti i complessi $((F_iC/F_{i+1}C)^*, d_i^*)$, ovvero che
%la mappa
%$$
%{\pi_i}^*_+: H^*((F_iC/F_{i+1}C)^* \otimes_R M) \to H^*((F_iC/F_{i+1}C)^* \otimes_R M_+)
%$$
%sia un isomorfismo.

If $C^* = \{0\}$
% o $C^* = R$
the assertion is obvious.
Suppose that $C^*$ has $m$ generators, with $m > 1$. Then the complexes $((F_iC/F_{i+1}C)^*, d_i^*)$
have a smaller number of generators and for $i \neq n-1, n$ they are well filtered. Therefore
we can suppose by induction that the map ${\pi_i}_+$, defined analogously to $\pi_+$,
induces an isomorphism in cohomology for all the complexes $((F_iC/F_{i+1}C)^*, d_i^*)$, $i \neq n-1, n$,
that is the map
$$
{\pi_i}^*_+: H^*((F_iC/F_{i+1}C)^* \otimes_R M) \to H^*((F_iC/F_{i+1}C)^* \otimes_R M'_+)
$$
is an isomorphism for such $i$.

The filtration $F$ on $C^*$ induces filtrations on $C_M^*$ and $C'^*_+$ in the following way:
$F_iC_M = F_iC \otimes_R M$, $F_iC'_+ = F_iC \otimes_R M'_+$. We have the following natural
isomorphisms:
$$
(F_iC/F_{i+1}C)^* \otimes_R M \simeq (F_iC_M/F_{i+1}C_M)^*
$$
$$
(F_iC/F_{i+1}C)^* \otimes_R M'_+ \simeq (F_iC'_+/F_{i+1}C'_+)^*.
$$
Through these isomorphisms the maps
$$(F_iC_M/F_{i+1}C_M)^* \to (F_iC'_+/F_{i+1}C'_+)^*$$
 induced by $\pi_+$
correspond to ${\pi_i}_+$ and hence induce an isomorphism in cohomology for $i \neq n-1, n$.

Let us consider the spectral sequences $E^{i,j}_r$ and $\overline{E}^{i,j}_r$ associated to
the complexes $C^*_M$ and ${C'}^*_+$ with the respective filtrations. We write $\pi_+^*$
also for the spectral sequences homomorphism induced by $\pi_+$. By the definition of 
the filtration $F$ we have that $E^{i,j}_r = \overline{E}^{i,j}_r = 0$ if $i >n$ or if $i=n, n-1$ and $j \neq 0$.
It is also clear that $E^{n-1,0}_1 \simeq E^{n,0}_1 = M$ and
$\overline{E}^{n,0}_1 \simeq \overline{E}^{n-1,0}_1 = M'_+$. For $0 \leq i<n-1$ we get that
$E^{i,j}_1 \simeq H^{i+j}(F_iC^*_M/F_{i+1}C^*_M)$ and
$\overline{E}^{i,j}_1 \simeq H^{i+j}(F_i{C'}^*_+/F_{i+1}{C'}^*_+)$ therefore the inductive hypothesis
gives that $E^{i,j}_1 \simeq \overline{E}^{i,j}_1$ and the isomorphism between the
terms of the spectral sequences is given by $\pi_+^*$.
Now consider the maps $d_1^{n-1,0}:M \to M$ and $\overline{d}_1^{n-1,0}: M'_+ \to M'_+$. By
%the inductive hypothesis on the complex $C^*$ 
condition (c)
we have that these maps correspond to the multiplication by a non-zero polynomial
$p = \sum_{i=s}^t b_i q^i $ with $b_s, b_t$ invertible elements of the ring $A$.
We can rewrite $p$ as follows:
$$
p = q^s b_s (1 + q p') = q^t b_t (1 + q^{-1} p'')
$$
with $p' \in A[q]$, $p'' \in A[q^{-1}]$.
%with $p', p'' \in A[x]$.
Now we can look at these elements in $M$:
$$
p^{-1}_+ = q^{-s} b_s^{-1} \sum_{i=0}^{\infty}(-qp')^i
$$
$$
p^{-1}_- = q^{-t} b_t^{-1} \sum_{i=0}^{\infty}(-q^{-1}p'')^i.
$$
Let $m \in M$, $m = \sum_{i \in \Z} a_i q^i$, we can write $m = m_+ + m_-$,
with $m_+ = \sum_{i=0}^{\infty} a_i q^i$ and $m_- = m - m_+$.
Notice that the products $p^{-1}_+ m_+$ and $p^{-1}_- m_-$ are well defined and the following equality holds:
$$
m = p (p^{-1}_+ m_+ + p^{-1}_- m_-).
$$
It turns out that the map $d_1^{n-1,0}:M \to M$ is surjective and the same holds, when passing to the
quotient, for the map $\overline{d}_1^{n-1,0}: M_+ \to M_+$.

Let us suppose that an element $m = \sum_{i \in \Z} a_i q^i$ is in the kernel of $d_1^{n-1,0}$.
This means that $pm=0$, that is for all integers $k$ we have:
$$
\sum_{i=s}^t b_i a_{k-i} = 0
$$
and so we obtain:
\begin{eqnarray}
a_k = -b_s^{-1} \sum_{i=1}^{t-s} b_{s+i} a_{k-i} \label{eqn1}\\
a_k = -b_t^{-1} \sum_{i=1}^{t-s} b_{t-i} a_{k+i}.\label{eqn2}
\end{eqnarray}
Therefore if we know a sequence of $t-s$ consecutive coefficients of an element $m$
sent to zero by the multiplication by $p$ we can use (\ref{eqn1}) and (\ref{eqn2}) to calculate recursively
all the other coefficients, determining $m$ completely.
So we find a bijection between $\ker d_1^{n-1,0}$ and
$\ker \overline{d}_1^{n-1,0}$. In fact, if $m \in M$ is such that $pm = 0$, then trivially also
$p[m]_+ = 0$ (we write $[m]_+$ for the equivalence class of $m$ in $M'_+$).
Conversely if $p[m]_+ = 0$ then
%
%correzione, ampliamento
%
we have $pm = z$, with $z \in A[q][[q^{-1}]]$, that is $z = \sum_{i \in \Z}v_iq^i$ with $v_i \in A$ and there exists an integer $l$
% \in \Z$ 
such that $v_i = 0$ for all $i > l$.
We can define recursively, for $j \geq 0$, the following elements:

$$
{\widetilde{a}[-1]}_i = a_i
$$
$$
{\widetilde{a}[j]}_i = \left\{
\begin{array}{lcr}
{\widetilde{a}[j-1]}_i & \mbox{\qquad if } & i \neq l - t -j \\
-b_t^{-1} \sum_{k=1}^{t-s} b_{t-k} {\widetilde{a}[j-1]}_{i+k} & \mbox{\qquad if } & i = l - t -j
\end{array}  \right.
$$
and
$$
{\widetilde{a}}_i = \left\{
\begin{array}{lcr}
a_i & \mbox{\qquad if } & i > l - t \\
{\widetilde{a}[l-t-i]}_{i} & \mbox{\qquad if } & i \leq l - t
\end{array}  \right.
$$
Notice that the coefficients $v_i$ for $i > h$ depend only on the coefficients $a_i$ for $i > h-t$, so if we write $\widetilde{m} = \sum_{i \in \Z}{\widetilde{a}}_i q^i$ we have that $p\widetilde{m} = 0$ and $[m]_+ = [\widetilde{m}]_+$.
%
%there exists one and only one element $\widetilde{m} \in M$ such that
%$p\widetilde{m} = 0$ and $[m]_+ = [\widetilde{m}]_+$.
%
%

To sum up we have that the map $\pi_+^*$ gives an isomorphism between the terms $E^{i,j}_1$ and
$\overline{E}^{i,j}_1$ for $i < n-1$ and between $\ker d_1^{n-1,0}$ and  $\ker \overline{d}_1^{n-1,0}$.
Moreover $E^{i,j}_2 = \overline{E}^{i,j}_2 = 0$ for $i = n-1$ and $j \neq 0$ and for $i > n-1$;
$\pi_+^*$ commutes with the differentials in the spectral sequences
(i. e. $\pi_+^* d_i = \overline{d}_i \pi_+^*$). We remark that
$\im d^{n-2,0}_1 \subset \ker d_1^{n-1,0}$ and
$\im \overline{d}^{n-2,0}_1 \subset \ker \overline{d}_1^{n-1,0}$ and so $\pi_+^*$ induces
an isomorphism between $\im d^{n-2,0}_1$ and $\im \overline{d}^{n-2,0}_1$. This implies that 
$\pi_+^*$ gives the isomorphisms
$E_2^{n-2,0}\simeq \overline{E}^{n-2,0}_2$ and $E_2^{n-1,0}\simeq \overline{E}^{n-1,0}_2$.
% and between $E_2^{n-1,0}$ and $\overline{E}^{n-1,0}_2$. 
Then we have a complete isomorphism between $E_2$ and $\overline{E}_2$ and
so between $E_\infty$ and $\overline{E}_\infty$. It follows that $\pi_+^*$ induces an
isomorphism in cohomology.

It is clear that the same fact holds for the map
% sarebbe potuto evitare di usare l'ipotesi che il
% coefficiente $b_t$ sia invertibile. Questa ?necessaria per ripetere il ragionamento,
$\pi_-: C^*_M \to C'^*_-$ and so Lemma is proved. \qed
\end{dm}
We write $\Phi$ for the isomorphism built in the proof of the previous Lemma.

\begin{lem} \label{l:diag}
In the exact sequence $(*)$ the map $\pi^*$ composed with the isomorphism $\Phi$
%del lemma \ref{l:iso}
corresponds to the diagonal map $\diag$:
$$
H^i(C^*_M) \stackrel{\diag}{\into} H^i(C^*_M) \oplus H^i(C^*_M).
$$
\end{lem}

\begin{dm}
It is enough to notice that, making the identification $H^*({C'}^*) = H^*({C'}^*_+) \oplus H^*({C'}^*_-)$, we have that
$\pi^* = (\pi_+^*,\pi_-^*)$ and so the statement follows immediately. \qed
\end{dm}

\begin{dm}[Proof (of Theorem \ref{t:centrale}).]
First of all we notice that, being $\pi^*$ injective,
$i^*$ turns out to be the null map and $\delta^*$ is surjective.
We call $p_1: H^i(C^*_M) \oplus H^i(C^*_M) \to H^i(C^*_M)$ the projection on the first component, $p_2$ the projection
on the second component and
$i_1: H^i(C^*_M) \into H^i(C^*_M) \oplus H^i(C^*_M)$ the inclusion defined by $i_1:b \mapsto (b,0)$.
Finally we define $\alpha = \delta^* \circ \Phi^{-1} \circ i_1$. We have the following diagram:
\begin{center}
\begin{tabular}{c}
\xymatrix
@R=1.5pc @C=2pc
{
0 \ar[r] & H^i(C^*_M) \imm[r]^-{\diag} \ar@{=}[d] &
 H^i(C^*_M) \oplus H^i(C^*_M) \ar@<0.5ex>[r]^-{p_1-p_2} \eq[d] \ar@{}[d]^\Phi &
 H^i(C^*_M) \ar@<0.5ex>[l]^-{i_1} \ar[r] \ar[d]^{\alpha} & 0 \\
0 \ar[r] & H^i(C^*_M) \ar[r]^-{\pi^*} & H^i({C'}^*_M) \ar[r]^-{\delta^*} & H^{i+1}(C^*) \ar[r] & 0 \\
}
\end{tabular}
\end{center}
Clearly both the lines are exact. We want to show that the diagram commutes.
The commutativity for the first square follows by Lemma \ref{l:diag},
so it remains to prove that the second square commutes.
A pair $(a,b) \in H^i(C^*_M) \oplus H^i(C^*_M)$ is sent, by the multiplication by $p_1 - p_2$, into
the element $a-b \in H^i(C^*_M)$. Then we have $i_1(a-b) = (a-b, 0)$ and the difference $(a,b) - (a-b,0) = (b,b)$ is
in the image of the map $\diag$. Therefore, because of the commutativity of the first square,
the images of the pairs $(a,b)$ and of $(a-b,0)$ in $H^i({C'}^*_M)$ are taken into the same
element by the map $\delta^*$. So we get the commutativity of the diagram.
The Theorem follows from the five lemma. \qed
\end{dm}

\section{Applications} \label{s:appl}
Let us consider a finite set $\Gamma$ endowed with a fixed total ordering.
We will indicate by $\Delta$ a generic subset of $\Gamma$. We also set again $R=A[q,q^{-1}]$,
with $A$ a commutative unitary ring. For every pair $(\Delta, w)$
with $\Delta \subset \Gamma$, $w \in \Gamma \setminus \Delta$ we associate a polynomial
$p_{\Delta, w}(q,q^{-1}) \in R \setminus \{ 0 \}$ such that the first and the last non-zero coefficients are
invertible in $A$.
Let also suppose that for every pair $(w,w')$ with $w \neq w'$ and $w,w' \in \Gamma \setminus \Delta$
the following equation holds:
\begin{equation}
p_{\Delta,w}(q,q^{-1}) p_{\Delta \cup \{w\},w'}(q,q^{-1}) + p_{\Delta,w'}(q,q^{-1}) p_{\Delta \cup \{w'\},w}(q,q^{-1}) = 0 \label{e:cobordo}
\end{equation}
Then we can consider the complex $(C_\Gamma^*,d^{*})$ defined as follows:
$$
C^*_\Gamma = \bigoplus_{\Delta \subset \Gamma} R.e_\Delta
$$
$$
d^*e_\Delta = \sum_{w \in \Gamma \setminus \Delta} p_{\Delta,w}(q,q^{-1}) e_{\Delta \cup \{w\}}.
$$
We remark that the relation (\ref{e:cobordo}) gives ${d^*}^2 = 0$.
We can also give a natural graduation to $C^*_\Gamma$ by defining the degree of an element
$e_\Delta$ as the cardinality of $\Delta$, so we get a cochain complex.

%Let $n$ be the cardinality of $\Gamma$. In the following we use the convention to 
%identify an element $e_\Delta \in C^*_\Gamma$ with the corresponding 
%$n$-string made up with $0$'s and $1$'s; here the $1$'s occupy the position 
%of the elements of $\Delta$ with respect to the given ordering of $\Gamma$. 
%So 
Without loss of generality we can think $\Gamma = \{1, \ldots, n\}$. 
We introduce the following notation: indicate by $\Gamma_i$ and $\Delta_i$ respectively the
subsets $\{ 1, \ldots, n-i-1\}$ and $\{ n-i+1, \ldots, n\}$.
We can filter the complex $C^*_\Gamma$ in the following way (see also
\cite{dps}): let $F_iC_\Gamma$ be the subcomplex generated by the elements $e_\Delta$, with $\Delta_i \sst \Delta$.

We have the following result:
\begin{teo} \label{t:appl}
With the filtration defined above the complex $(C^*_\Gamma, d^*)$ is well filtered.
\end{teo}
\begin{dm}
We can prove this by induction on the cardinality of $\Gamma$. If $\Gamma$ is empty the Theorem is obvious.
Therefore let us suppose that the Theorem holds for all the complexes made up from a set with
less than $n$ elements and we prove it for a complex $C^*_\Gamma$, with $\Gamma = \{1, \ldots, n\}$.

It is straightforward to see that $F_0C_\Gamma = C^*_\Gamma$ and $F_{n+1}C_\Gamma = \{ 0 \}$. Moreover
$F_nC_\Gamma$ and $F_{n-1}C_\Gamma / F_nC_\Gamma$ are generated respectively by the elements $e_\Gamma$ and
$e_{\Delta_{n-1}}$ and they are both isomorphic to $R$. The induced differential
$$
\overline{d}: F_{n-1}C_\Gamma/F_nC_\Gamma \to F_nC_\Gamma/F_{n+1}C_\Gamma
$$
corresponds to the multiplication by the polynomial $p_{\Delta_{n-1},1}(q,q^{-1})$.

Finally the complex $((F_iC_\Gamma/F_{i+1}C_\Gamma)^*, d_i^*)$ is isomorphic to the complex
$C^*_{\Gamma_i}$, where the coboundary is defined by the polynomials
$$
\overline{p}_{\Delta, j}(q,q^{-1}) := p_{\Delta \cup \Delta_i,j}(q,q^{-1}) \qquad \mbox{for }\Delta \subset \Gamma_i,
j \in \Gamma_i \setminus \Delta
$$
and so it is well filtered by induction. \qed
\end{dm}

Now we apply last result and Theorem \ref{t:centrale} to the cohomology with local
coefficients of Artin groups. In \cite{salv1} Salvetti proved that:

\begin{teo} \label{t:salv}
Let $W$ be a Coxeter group with generating set $\Gamma$ with a fixed total ordering and let $G_W$ be the associated Artin group. Let $R$ be a commutative ring
with unit and let $q$ be a unit in $R$ and let $M$ be an $R$-module. We write $W_{\Delta}(q)$ for the
Poincar\'e polynomial of the subgroup of $W$ generated by $\Delta$, with $\Delta \sst \Gamma$.
Let $\cl{L}_q = \cl{L}_q(X_W; M)$ be the local system over
$G_W$ with coefficients in $M$ given by the map that sends every standard generator
of $G_W$ into the automorphism of $M$ given by the multiplication by $q$. Then
$$
H^*(G_W; \cl{L}_q) \simeq H^*(C^*)
$$
where
$$
C^k= \{ \sum a_{\Delta} e_\Delta \mid a_{\Delta} \in M, \Delta \sst \Gamma, | \Delta | = k \}
$$
and the coboundary is given by
$$
\delta^k(e_\Delta) = \sum_{j \in \Gamma \setminus \Delta} {(-1)}^{\sigma(j,\Delta)}
\frac{W_{\Delta \cup \{j\}}(-q)}{W_{\Delta}(-q)}e_{\Delta \cup \{j\}}
$$
where $\sigma(j,\Delta) = | \{ i \in \Delta, i < j \} |$. \qed
\end{teo}

%\begin{os}
\begin{prop}
Let $R = A[q,q^{-1}]$ and $M = R$. Then the complex $C^*$ in Theorem \ref{t:salv} is well filtered.
\end{prop}
\begin{dm}
In fact the polynomial $W_\Gamma(q)$ divides $W_{\Gamma'}(q)$ when $\Gamma \subset \Gamma'$.
Moreover the polynomials
$W_\Gamma(q)$ are products of cyclotomic polynomials (see \cite{bourb}),
so they have first and last non-zero coefficients
equal to $1$. By using Theorem \ref{t:appl} we can easily see that 
%for $R = A[q,q^{-1}]$ and $M = R$ the complex
$C^*$ 
%in Theorem \ref{t:salv} 
is well filtered. \qed
\end{dm}
%As a consequence we get that
%$$
%H^*(G_W,\cl{L}_q') = H^{*+1}(G_W,\cl{L}_q). 
%H^*(G_W,A[[q,q^{-1}]]) = H^{*+1}(G_W,A[q,q^{-1}]). 
%$$
%\end{os}

Now let $W$ be a finite Coxeter group. We can think of $W$ as a group generated by orthogonal reflections
in a real vector space $\B{V}$. Let $\g{H}$ be the arrangement of all the hyperplanes in $\B{V}$ such that
the associated orthogonal reflection is in $W$. We can consider the complexified space $\B{V}_\C$ and
the complexified arrangement $\g{H}_\C$. For every hyperplane $\B{H} \in \g{H}_\C$ we chose a linear
function $l_{\B{H}}$ such that $\ker l_{\B{H}} = \B{H}$. The polynomial
$$\delta = \prod_{\B{H} \in \g{H}_\C} l_{\B{H}}^2 $$
is called the discriminant of the arrangement and it is invariant with respect to the diagonal action of $W$ on $\B{V}_\C$.
The space $$\B{X}_W = (\B{V}_\C \setminus \cup_{\B{H} \in \g{H}}\B{H})/W $$ is a classifying space
for the Artin group $G_W$ (see \cite{deligne}), and $\delta$ induces a fibering $$\delta':\B{X}_W \to \C^*.$$
The fiber $\B{F}_W = {\delta}^{-1}(1)$ is called the Milnor fiber of $\B{D}_W = (\cup_{\B{H} \in \g{H}}\B{H})/W$.
The associated homotopy exact sequence gives us
that the $\B{F}_W$
is a classifying space for the subgroup $H_W < G_W$, which is the kernel of the natural homomorphism
$$G_W \to \Z $$
defined by sending each standard generator to $+1$.

Now we set again $R = A[q,q^{-1}]$ and $C^*$ and let $C^*_M = C^* \otimes M$
be the algebraic complexes defined as in Theorem \ref{t:salv}, over $R$ or $M = A[[q,q^{-1}]]$ respectively.
Then (by definition) the following equality holds:
\begin{equation} \label{e:prima}
H^*(\B{F}_W;A) = H^*(H_W;A)
\end{equation}
and the Shapiro Lemma (see \cite{brown}) gives that
\begin{equation} \label{e:seconda}
H^*(H_W;A) = H^*(G_W;Coind_{H_W}^{G_W}A) = H^*(G_W;M) = H^*(C^*_M)
\end{equation}
where the action of $G_W$ over $M$ is given by sending each standard generator into
the multiplication by $q$.
From Theorem \ref{t:centrale} and the remark following Theorem \ref{t:salv}, we get that
\begin{equation} \label{e:terza}
H^*(G_W,M) = H^{*+1}(G_W,R)
\end{equation}

Using equalities (\ref{e:prima}), (\ref{e:seconda}) and (\ref{e:terza})
%what we have seen in theorems \ref{t:centrale}, \ref{t:appl} and \ref{t:salv}
we get immediately the following result:

\begin{teo} \label{t:finale}
Let $W$ be a finite irreducible Coxeter group and let
$$
\B{F}_W \into \B{X}_W \stackrel{\delta'}{\to} \C^*
$$
be the fibration defined as before. Let $R = A[q,q^{-1}]$ be considered as a
$G_W$-module with the action defined before. Then the following equality holds:
$$
H^*(\B{F}_W;A) = H^{*+1}(G_W; R)
$$
\qed
\end{teo}

\section*{Acknowledgment}
It is a pleasure to thank Mario Salvetti for his collaboration and his help.

\end{document}